\documentclass[a4,10pt]{article}
\usepackage{amsmath,amscd,amssymb}

\newcommand{\nbigh}{\mathcal{H}}

\newcommand{\nbign}{\mathcal{N}}
\newcommand{\nbigo}{\mathcal{O}}

\newcommand{\nbigu}{\mathcal{U}}
\newcommand{\nbigv}{\mathcal{V}}

\newcommand{\cnum}{{\mathbb C}}
\newcommand{\real}{{\mathbb R}}

\newcommand{\vece}{{\boldsymbol e}}

\newcommand{\lrarr}{\longrightarrow}

\newcommand{\pf}{{\bf Proof}\hspace{.1in}}
\newcommand{\qed}{\mbox{\rule{1.2mm}{3mm}}}

\def\Re{\mathop{\rm Re}\nolimits}

\def\rank{\mathop{\rm rank}\nolimits}

\def\Ker{\mathop{\rm Ker}\nolimits}

\def\Tr{\mathop{\rm Tr}\nolimits}

\def\id{\mathop{\rm id}\nolimits}

\newcommand{\del}{\partial}
\newcommand{\delbar}{\overline{\del}}

\newcommand{\Etilde}{\widetilde{E}}

\newcommand{\nablatilde}{\widetilde{\nabla}}

\newcommand{\wbar}{\overline{w}}

\newcommand{\htilde}{\widetilde{h}}

\newcommand{\vtilde}{\widetilde{v}}

\newcommand{\stilde}{\widetilde{s}}

\newcommand{\Utilde}{\widetilde{U}}

\newcommand{\ubar}{\overline{u}}

\newcommand{\etilde}{\widetilde{e}}

\newcommand{\vecsigma}{\boldsymbol\sigma}

\newtheorem{thm}{Theorem}[section]

\newtheorem{rem}[thm]{Remark}
\newtheorem{lem}[thm]{Lemma}
\newtheorem{prop}[thm]{Proposition}
\newtheorem{df}[thm]{Definition}

\newtheorem{example}[thm]{Example}

\begin{document}

\title{Some characterizations of Dirac type singularity
 of monopoles}

\author{Takuro Mochizuki and Masaki Yoshino}

\date{}
\maketitle

\begin{abstract}
We study singular monopoles
on open subsets in the $3$-dimensional Euclidean space.
We give two characterizations of Dirac type singularities.
One is given in terms of the growth order of the norms of sections
which are invariant by the scattering map.
The other is given in terms of the growth order of the norms of the Higgs
fields.

\vspace{.1in}

\end{abstract}

\section{Introduction}

\subsection{Dirac type singularity}
Let $U$ be any open subset in $\real\times\cnum$
such that $(0,0)\in U$.
We regard $U$ as the Riemannian manifold
with the standard Euclidean metric
$dt\,dt+dw\,d\wbar$,
where $(t,w)$ is the standard coordinate of $\real\times\cnum$.
We set $U^{\ast}:=U\setminus\{(0,0)\}$.
Let $(E,h,\nabla,\phi)$ be a monopole on $U^{\ast}$.
Namely, 
$E$ is a $C^{\infty}$-bundle on $U^{\ast}$
with a Hermitian metric $h$,
a unitary connection $\nabla$,
and an anti-self-adjoint endomorphism $\phi$
satisfying the Bogomolny equation
$F(\nabla)=\ast \nabla\phi$.
Here, $F(\nabla)$ is the curvature of $\nabla$,
and $\ast$ is the Hodge star operator.
We regard $(E,h,\nabla,\phi)$
as a singular monopole on $U$.

The study of singular monopole was pioneered 
by Kronheimer  \cite{Kronheimer-Master-Thesis}.
Among other things,
he introduced a reasonable class of singularity,
called Dirac type singularity.
Let $\varphi:\cnum^2\lrarr\real\times\cnum$
be given by
$\varphi(u_1,u_2)=(|u_1|^2-|u_2|^2,2u_1u_2)$.
We set
$(\Etilde,\htilde):=
\varphi^{\ast}(E,h)$
and 
$\nablatilde:=
 \varphi^{\ast}\nabla
+\sqrt{-1}\xi\otimes\varphi^{\ast}\phi$
on $\varphi^{-1}(U^{\ast})$,
where 
$\xi:=-u_1d\ubar_1+\ubar_1du_1-\ubar_2du_2+u_2d\ubar_2$.
As discovered by Kronheimer,
$(\Etilde,\htilde,\nablatilde)$
is an instanton on $\varphi^{-1}(U^{\ast})$,
i.e.,
the curvature $F(\nablatilde)$ is a $(1,1)$-form
and satisfies $\Lambda_{\cnum^2}F(\nablatilde)=0$.
Then, the point $(0,0)\in U$ is called a Dirac type singularity of
the monopole $(E,h,\nabla,\phi)$
if $(\Etilde,\htilde,\nablatilde)$ is extended to an instanton
on $\varphi^{-1}(U)$.

The condition restricts the behaviour of 
the monopole around $(0,0)$.
Set $R(t,w):=\sqrt{|t|^2+|w|^2}$.
In the case of $SU(2)$-monopoles,
according to Kronheimer \cite{Kronheimer-Master-Thesis},
if $(0,0)$ is a Dirac type singularity,
the limit $\lim_{(t,w)\to (0,0)}\bigl|(R\phi)_{(t,w)}\bigr|$
exists, and $\nabla(R\phi)$ is bounded.
Moreover, he proved the converse,
i.e., the conditions imply
that the point $(0,0)$ is a Dirac type singularity of the monopole.
See \cite{Pauly} for a generalization
to the context of general Riemannian three manifolds.
See \cite{Charbonneau-Hurtubise}
for the higher rank case.
In {\rm\cite{Cherkis-Kapustin2}},
Cherkis and Kapustin gave 
a characterization of Dirac type singularity
in terms of the growth order of
the Higgs field and the curvature,
i.e.,
 $\bigl|\phi\bigr|_h=O(R^{-1})$
and 
 $\bigl|F(\nabla)\bigr|_h=O(R^{-2})$.
It is particularly useful 
in the study of Nahm transforms.

\subsection{Main results}
\label{subsection;17.2.20.1}

In this paper, we study two other characterizations 
of Dirac type singularity.
The first one is given in terms of the growth order of 
the norms of sections which are invariant by scattering map.
The second one is given in terms of the growth order of 
the norms of the Higgs field.
The latter can be stated in a simple way
as follows.

\begin{thm}[Theorem
 \ref{thm;16.12.14.10}]
 \label{thm;16.12.19.2}
The point $(0,0)$ is a Dirac type singularity of 
the monopole
$(E,h,\nabla,\phi)$
if and only if 
$|\phi|_h=O(R^{-1})$.
\end{thm}

To state the former characterization,
we give preliminaries.
For simplicity, we suppose that
$U$ is the product of
connected open subsets
$U_t\subset\real$ and $U_w\subset\cnum$.
Set $U_{+}:=\bigl\{(t,w)\in U\,\big|\,t\geq 0\bigr\}$
and $U_{-}:=\bigl\{(t,w)\in U\,\big|\,t\leq 0\bigr\}$.

Take a small $\epsilon>0$
such that $\{\pm\epsilon\}\subset U_t$.
We put
$E^{\epsilon}:=E_{|\{\epsilon\}\times U_w}$
and 
$E^{-\epsilon}:=E_{|\{-\epsilon\}\times U_w}$.
Then, for any $C^{\infty}$-section $s^{\pm\epsilon}$
of $E^{\pm\epsilon}$,
we have a unique $C^{\infty}$-section 
$\stilde^{\pm\epsilon}$ of $E_{|U_{\pm}\setminus\{(0,0)\}}$
satisfying
 $(\nabla_t-\sqrt{-1}\phi)\stilde^{\pm\epsilon}=0$
and
$\stilde^{\pm\epsilon}_{|\{\pm\epsilon\}\times U_w}
=s^{\pm\epsilon}$.

The dual bundle $E^{\lor}$ of $E$
is naturally equipped with 
the induced metric $h_{E^{\lor}}$,
the induced unitary connection $\nabla_{E^{\lor}}$.
Let $\phi^{\lor}$ denote the endomorphism of $E^{\lor}$
obtained as the dual of $\phi$.
Then,
$(E^{\lor},h_{E^{\lor}},\nabla_{E^{\lor}},-\phi^{\lor})$
is also a monopole on $U^{\ast}$.

We say that 
$(E,h,\nabla,\phi)$ satisfies the condition {\bf (D)}
if the following holds:
\begin{itemize}
\item
For any $C^{\infty}$-section $s^{\pm\epsilon}$
of $E^{\pm\epsilon}$,
we have 
$\bigl|
 \stilde^{\pm\epsilon}
 \bigr|_h=O\bigl((|t|+|w|)^{-N}\bigr)$
on $U_{\pm}\setminus\{(0,0)\}$
for some $N>0$.
\item
The same estimate holds for sections of
$(E^{\lor},h_{E^{\lor}},\nabla_{E^{\lor}},-\phi^{\lor})$.
\end{itemize}
Then, the former characterization is stated as follows.
\begin{thm}[Theorem 
 \ref{thm;13.12.10.10}]
\label{thm;16.12.19.1}
The point $(0,0)$ is a Dirac type singularity
of the monopole $(E,h,\nabla,\phi)$
if and only if the condition {\bf (D)} is satisfied.
\end{thm}

For the proof of Theorem \ref{thm;16.12.19.1},
we apply a deep result of Donaldson 
on the Dirichlet problem for instantons
\cite[Theorem 1]{Donaldson-boundary-value}.
Theorem \ref{thm;16.12.19.2}
is an easy consequence of 
Theorem \ref{thm;16.12.19.1}.

\vspace{.1in}

We are motivated by the study of the Nahm transform,
which produces singular instantons on $\real^4$
equivariant with respect to 
the actions of a subgroup of $\real^4$,
from other type of equivariant instantons.
Here, we naturally regard monopoles, harmonic bundles,
and solutions of the Nahm equation as equivariant instantons.
For instance, 
it produces singular monopoles on $S^1\times\real^2$
from harmonic bundles on $S^1\times\real$
\cite{Cherkis-Kapustin1,Cherkis-Kapustin2},
and singular monopoles on $(S^1)^3$
from instantons on $(S^1)^3\times\real$
\cite{Charbonneau}.
(More precisely, we need to impose boundary conditions
on corresponding harmonic bundles and instantons,
but we omit such details here.)
To prove that the singularities of the monopoles
are of Dirac type,
it is convenient to have characterizations 
which can be checked rather easily.
Indeed, Cherkis and Kapustin \cite{Cherkis-Kapustin2}
gave their characterization for that purpose.
The first author obtained Theorem \ref{thm;16.12.19.1}
in the study of the Nahm transform from
singular harmonic bundles on $S^1\times\real$
to singular monopoles on $S^1\times\real^2$
\cite{Mochizuki-periodic-monopole},
and the authors obtained Theorem \ref{thm;16.12.19.2}
in the study of the Nahm transform from
instantons on $(S^1)^3\times\real$
to singular monopoles on $(S^1)^3$ \cite{Yoshino}.

\section{Mini-holomorphic bundles and holomorphic bundles}

\subsection{Mini-holomorphic bundles}

Let $t$ and $w$ be the standard coordinate
of $\real$ and $\cnum$,
respectively.
We have the real vector field $\del_t$
and a complex vector field $\del_{\wbar}$
on $\real\times\cnum$.
A $C^{\infty}$-function $f$ on an open subset 
$U\subset\real\times\cnum$
is called mini-holomorphic
if $\del_tf=0$ and $\del_{\wbar}f=0$.
(We use the  prefix ``mini'' 
by following ``mini-twistor'' 
in \cite{Kronheimer-Master-Thesis}.)

Let $E$ be a $C^{\infty}$-vector bundle
on an open subset $U\subset\real\times\cnum$.
Let $C^{\infty}(E)$ denote the space of
$C^{\infty}$-section of $E$.
A mini-holomorphic structure of $E$
is a pair of differential operators
$(\del_{E,t},\del_{E,\wbar})$
on $C^{\infty}(E)$
satisfying the following conditions:
\begin{itemize}
\item
 For any $f\in C^{\infty}(U)$ and $s\in C^{\infty}(E)$,
 we have
 $\del_{E,t}(fs)=\del_t(f)\cdot s+f\,\del_{E,t}(s)$
 and 
 $\del_{E,\wbar}(fs)=\del_{\wbar}(f)\cdot s+f\,\del_{E,\wbar}(s)$.
\item
 We have the commutativity
 $[\del_{E,t},\del_{E,\wbar}]=0$.
\end{itemize}
A mini-holomorphic bundle
means a $C^{\infty}$-bundle
with a mini-holomorphic structure
$(E,\del_{E,t},\del_{E,\wbar})$.
A $C^{\infty}$-section $s$ of $E$ is called mini-holomorphic
if $\del_{E,t}s=0$ and $\del_{E,\wbar}s=0$.

\begin{rem}
The concept of mini-holomorphic bundles
was efficiently used 
in previous studies
{\rm\cite{Charbonneau-Hurtubise}},
{\rm \cite{Cherkis-Kapustin1,Cherkis-Kapustin2}},
{\rm \cite{Hitchin-monopoles-geodesics,Hitchin-construction-monopole}},
{\rm \cite{Kronheimer-Master-Thesis}}
and {\rm\cite{Norbury}}, etc.
\hfill\qed
\end{rem}

Suppose that $U$ is the product of 
an interval $U_t\subset\real$
and an open subset $U_w\subset\cnum$.
Let $(E,\del_{E,t},\del_{E,\wbar})$
be a mini-holomorphic bundle on $U$.
For each $t\in U_t$,
we set 
$E^{t}:=E_{|\{t\}\times U_w}$,
which is equipped with a naturally induced
holomorphic structure $\delbar_{E^{t}}$.
For any $w\in U_w$,
we have the naturally induced connection
of $E_{|U_t\times\{w\}}$.
We have the parallel transport
$\Phi_w^{t_2,t_2}:E_{(t_1,w)}\lrarr E_{(t_2,w)}$.
They give an isomorphism
of holomorphic bundles
$\Phi^{t_2,t_1}:
 (E^{t_1},\delbar_{E^{t_1}})
 \simeq
 (E^{t_2},\delbar_{E^{t_2}})$
on $U_w$.
It is called the scattering map
in \cite{Charbonneau-Hurtubise}.
For any $t\in U_t$,
the restriction induces the bijection
between
the mini-holomorphic sections of $E$
and the holomorphic sections of $E^t$.

\subsection{The induced holomorphic bundles}
\label{subsection;16.12.16.10}

Following \cite{Charbonneau-Hurtubise}
and \cite{Kronheimer-Master-Thesis},
we consider the map
$\varphi:\cnum^2\lrarr \real\times\cnum$
given by
$\varphi(u_1,u_2)
=\bigl(|u_1|^2-|u_2|^2,2u_1u_2\bigr)$.
Note that $\varphi^{-1}(0,0)=\{(0,0)\}$.
We have the $S^1$-action on $\cnum^2$
given by
$e^{\sqrt{-1}\theta}(u_1,u_2)
=(e^{\sqrt{-1}\theta}u_1,e^{-\sqrt{-1}\theta}u_2)$.
We can naturally identify 
$\varphi$
with 
the projection to the quotient space
$\cnum^2\lrarr \cnum^2/S^1$.

We have a naturally defined bundle map
$\varphi_{\ast}:
 T(\cnum^2\setminus\{(0,0)\})\otimes_{\real}\cnum
\lrarr 
 \varphi^{-1}T(\real\times\cnum)\otimes_{\real}\cnum$
induced by the tangent map.
We have the following formula
on $\cnum^2\setminus\{(0,0)\}$:
\[
 \varphi_{\ast}(\del_{\ubar_1})
=u_1\varphi^{-1}(\del_t)
+2\ubar_2\varphi^{-1}(\del_{\wbar}),
\quad\quad
 \varphi_{\ast}(\del_{\ubar_2})
=-u_2\varphi^{-1}(\del_t)
+2\ubar_1\varphi^{-1}(\del_{\wbar}).
\]
The following lemma is easy to see.
\begin{lem}
Let $U$ be an open subset in $\real\times\cnum$
such that $(0,0)\not\in U$.
A $C^{\infty}$-function $f$ on $U$
is mini-holomorphic
if and only if
$\varphi^{\ast}(f)$ is holomorphic on $\varphi^{-1}(U)$.
\hfill\qed
\end{lem}

Let $U$ be any open subset in $\real\times\cnum$
such that $(0,0)\not\in U$.
Let $(E,\del_{E,t},\del_{E,\wbar})$
be a mini-holomorphic bundle on $U$.
We set
$\Utilde:=\varphi^{-1}(U)\subset\cnum^2$.
We put $\Etilde:=\varphi^{-1}(E)$.
Let $C^{\infty}(\Etilde)$ denote the space of
$C^{\infty}$-sections of $\Etilde$ on $\Utilde$.

\begin{lem}
We have the unique differential operators
$\del_{\Etilde,\ubar_i}$ $(i=1,2)$
on $C^{\infty}(\Etilde)$
satisfying the following conditions.
\begin{itemize}
\item
 For any $s\in C^{\infty}(E)$
 and $f\in C^{\infty}(U)$,
 we have
\begin{equation}
 \label{eq;16.12.15.1}
  \del_{\Etilde,\ubar_1}(f\varphi^{-1}(s))
=\del_{\ubar_1}(f)\cdot \varphi^{-1}(s)
+u_1\varphi^{-1}(\del_{E,t}s)
+2\ubar_2\varphi^{-1}(\del_{E,\wbar}s),
\end{equation}
\begin{equation}
 \label{eq;16.12.15.2}
 \del_{\Etilde,\ubar_2}(f\varphi^{-1}(s))
=\del_{\ubar_2}(f)\cdot\varphi^{-1}(s)
-u_2\varphi^{-1}(\del_{E,t}s)
+2\ubar_1\varphi^{-1}(\del_{E,\wbar}s).
\end{equation}
\end{itemize}
Moreover,
we have the commutativity
$[\del_{\Etilde,\ubar_1},\del_{\Etilde,\ubar_2}]=0$.
\end{lem}
\pf
By the uniqueness,
it is enough to check the claim locally around any point of $U$.
Hence, we may assume to have a mini-holomorphic frame 
$v_1,\ldots,v_r$ of $E$.
We obtain the frame
$\varphi^{-1}(v_1),\ldots,\varphi^{-1}(v_r)$
of $\Etilde$.
We define the differential operators
$\del_{\Etilde,\ubar_i}$ as follows:
\begin{equation}
\label{eq;16.12.15.3}
 \del_{\Etilde,\ubar_i}\Bigl(\sum_{j=1}^r f_j\varphi^{-1}(v_j)\Bigr)
=\sum_{j=1}^r \del_{\ubar_i}(f)\cdot \varphi^{-1}(v_j).
\end{equation}
We can check 
the equalities (\ref{eq;16.12.15.1}), (\ref{eq;16.12.15.2}),
and the commutativity 
$[\del_{\Etilde,\ubar_1},\del_{\Etilde,\ubar_2}]=0$
easily. 
The uniqueness is clear.
\hfill\qed

\subsection{Extendability of the induced bundle}
\label{subsection;16.12.18.101}

Let $U_t\subset\real$ be a neighbourhood of $0$.
Let $U_w\subset\cnum$ be a neighbourhood of $0$.
We assume that $U_t$ and $U_w$ are connected.
We set $U:=U_t\times U_w$.
We set $U^{\ast}:=U\setminus\{(0,0)\}$.
We set
$A^-:=(U_t\times U_w)\setminus\{(t,0)\,|\,t\geq 0\}$
and 
$A^+:=(U_t\times U_w)\setminus\{(t,0)\,|\,t\leq 0\}$.

Let $(E,\del_{E,t},\del_{E,\wbar})$
be a mini-holomorphic bundle on $U^{\ast}$.
Take $\epsilon>0$
such that $\{\pm\epsilon\}\subset U_t$.
We have the scattering map
of holomorphic vector bundles
$\Phi^{\epsilon,-\epsilon}:
 E^{-\epsilon}_{|U_w\setminus\{0\}}
\simeq
 E^{\epsilon}_{|U_w\setminus\{0\}}$.

\begin{lem}
\label{lem;16.12.18.3}
The following conditions are equivalent.
\begin{description}
\item[(A)]
 $\Phi^{\epsilon,-\epsilon}$ is meromorphic at $0$.
 Namely, for any holomorphic section $s$ of $E^{-\epsilon}$,
 the induced section $\Phi^{\epsilon,-\epsilon}(s)$
 of $E^{\epsilon}_{|U_w\setminus\{0\}}$
 is a meromorphic section of
 $E^{\epsilon}$.
\item[(B)]
 The induced holomorphic bundle
 $\Etilde$ on $\varphi^{-1}(U^{\ast})$
 is extended to a holomorphic vector bundle
 on $\varphi^{-1}(U)$.
\end{description}
\end{lem}
\pf
Suppose the condition {\bf (A)}.
Let $\nbigv_{\pm}$
denote the sheaf of holomorphic sections of
$E^{\pm\epsilon}$.
Let $\nbigo_{U_w}(\ast 0)$
denote the sheaf of meromorphic functions
which admit poles at $0$.
We set
$\nbigu:=\nbigv_+\otimes\nbigo_{U_w}(\ast 0)
\simeq
 \nbigv_-\otimes\nbigo_{U_w}(\ast 0)$.
We set $\nbigv_0$
as the intersection of
$\nbigv_-$ and $\nbigv_+$
in $\nbigu$.
Set $r:=\rank E$.
By shrinking $U_w$,
we have a frame $e_1,\ldots,e_r$
of $\nbigv_0$,
non-negative integers
$m_{\pm,1},\ldots,m_{\pm,r}$
such that
(i) $w^{-m_{\pm,i}}e_i$ $(i=1,\ldots,r)$
give a frame of $\nbigv_{\pm}$,
(ii) $m_{+,i}m_{-,i}=0$ hold for any $i$.

We have the mini-holomorphic sections of $v_i$ of $E$ on $U^{\ast}$
corresponding to $e_i$ $(i=1,\ldots,r)$.
On $A^{\pm}$,
we have the mini-holomorphic sections
$w^{-m_{\pm,i}}v_i$ $(i=1,\ldots,r)$
which give a mini-holomorphic frame of 
$E_{|A^{\pm}}$.

We have the holomorphic sections
$\vtilde_i:=\varphi^{\ast}(v_i)$ of $\Etilde$
$(i=1,\ldots,r)$
on $\varphi^{-1}(U^{\ast})$,
which give a frame of $\Etilde$
on $\varphi^{-1}(U)\setminus
 \bigl(\{u_1u_2=0\}\bigr)$.
By the previous consideration,
around any point of 
$P\in \{(u_1,0)\}\cap\varphi^{-1}(U^{\ast})$,
we have the holomorphic frame
of $\Etilde$
given by 
$u_2^{-m_{+,i}}\vtilde_i$ $(i=1,\ldots,r)$.
Around any point
$P\in \{(0,u_2)\}\cap\varphi^{-1}(U^{\ast})$,
we have the holomorphic frame of $\Etilde$
given by 
$u_1^{-m_{-,i}}\vtilde_i$ $(i=1,\ldots,r)$.
Hence, the holomorphic bundle $\Etilde$ has 
a frame 
$u_1^{-m_{-,i}}u_2^{-m_{+,i}} \vtilde_i$
$(i=1,\ldots,r)$,
and $\Etilde$ is extended to a holomorphic bundle
on $\varphi^{-1}(U)$,
i.e., the condition {\bf (A)} is satisfied.

Suppose {\bf (B)}.
Let $(\Etilde_0,\delbar_{\Etilde_0})$
denote the holomorphic bundle
obtained as 
the extension of $(\Etilde,\delbar_{\Etilde})$
on $\varphi^{-1}(U)$.
We may assume that $U$ is bounded.
Note that we have the $S^1$-action
on $\varphi^{-1}(U)$,
and $(\Etilde,\delbar_{\Etilde})$ is $S^1$-equivariant.
Let us observe that
$(\Etilde_0,\delbar_{\Etilde_0})$
is naturally $S^1$-equivariant.
We take a holomorphic frame
$a_1,\ldots,a_r$
on an $S^1$-invariant neighbourhood
$\Utilde_0$ of $(0,0)$.
For any $b=e^{\sqrt{-1}\theta}\in S^1$,
we obtain a holomorphic frame
$b^{\ast}a_1,\ldots,b^{\ast}a_r$
of $\Etilde_{|\Utilde_0\setminus\{(0,0)\}}$.
By the Hartogs theorem,
$b^{\ast}a_j$ are holomorphic sections of
$\Etilde_{0|\Utilde_0}$.
We obtain the matrix valued functions
$A_{i,j}(b,u_1,u_2)$ determined by
$b^{\ast}a_j=\sum A_{i,j}(b,u_1,u_2)a_i$.
When $b$ is fixed,
$A_{i,j}(b,u_1,u_2)$ are holomorphic with respect to
$(u_1,u_2)$.
Because $\det(A_{i,j})$ is nowhere vanishing on
$\Utilde_0\setminus\{(0,0)\}$,
we obtain that $\det(A_{i,j})$
is nowhere vanishing on $\Utilde_0$,
i.e., $b^{\ast}a_1,\ldots,b^{\ast}a_r$
also give a frame of $\Etilde_{0|\Utilde_0}$.
Note that $A_{i,j}(b,u_1,u_2)$ are $C^{\infty}$
on $S^1\times (\Utilde_0\setminus\{(0,0)\})$.
By using that $A_{i,j}(b,u_1,u_2)$ are holomorphic
with respect to $(u_1,u_2)$,
we can easily deduce that
$A_{i,j}(b,u_1,u_2)$ are $C^{\infty}$
on $S^1\times \Utilde_0$.

We take an $S^1$-equivariant Hermitian metric
$h_0$ of $\Etilde_0$.
Let $\nbigh(\Etilde_0)$
denote the space of holomorphic sections
of $\Etilde_0$,
which is $L^2$
with respect to
$h_0$ and the volume form
associated to $du_1d\ubar_1+du_2d\ubar_2$.
The restriction map
$\Phi:\nbigh(\Etilde_0)
\lrarr
 \Etilde_{0|(0,0)}$
is surjective and $S^1$-equivariant.
The orthogonal complement of $\Ker\Phi$
is an $S^1$-subspace of 
$\nbigh(\Etilde_0)$,
and it is isomorphic to $\Etilde_{0|(0,0)}$
by the restriction.
We have a frame 
$s_1,\ldots,s_r$
of $\Etilde_{0|(0,0)}$
and integers $m_1,\ldots,m_r$
such that
$b^{\ast}(s_j)=b^{m_j}s_j$
for $b\in S^1$.
We have sections
$\stilde_j$ $(j=1,\ldots,r)$
such that
$b^{\ast}(\stilde_j)=b^{m_j}\stilde_j$.
By shrinking $U$,
we may assume that
$\stilde_1,\ldots,\stilde_r$
is a frame of $\Etilde$
on $\varphi^{-1}(U)$.

On $\varphi^{-1}(A^+)=\varphi^{-1}(U)\setminus\{u_1=0\}$,
we have 
the $S^1$-invariant frame
$u_1^{-m_j}\stilde_j$ $(j=1,\ldots,r)$,
which induces a holomorphic frame
$\sigma_{+,j}$ $(j=1,\ldots,r)$
on $E_{|A^+}$.
On $\varphi^{-1}(A^-)=\varphi^{-1}(U)\setminus\{u_2=0\}$,
we have 
the $S^1$-invariant frame
$u_2^{m_j}\stilde_j$ $(j=1,\ldots,r)$,
which induces a holomorphic frame
$\sigma_{-,j}$ $(j=1,\ldots,r)$
on $E_{|A^-}$.
The scattering map is given by
$\sigma_{-,j}
\longmapsto
 (w/2)^{m_j}\sigma_{+,j}$,
i.e.,
the condition {\bf (A)} is satisfied.
\hfill\qed

\section{Monopoles and instantons}

\subsection{Monopoles and underlying mini-holomorphic structure}

Let $U$ be any open subset in $\real\times\cnum$.
We regard $U$ as a Riemannian manifold
with the standard Euclidean metric.
Let $E$ be a $C^{\infty}$-bundle on $U$
with a Hermitian metric $h$ and a unitary connection $\nabla$.
Let $\phi$ be an anti-self-adjoint endomorphism of $E$.
Let $F(\nabla)$ denote the curvature of $\nabla$.
The tuple $(E,h,\nabla,\phi)$ is called a monopole
if the Bogomolny equation
$F(\nabla)=\ast \nabla\phi$ is satisfied,
where $\ast$ is the Hodge star operator.

We have the differential operator
$\del_{E,\wbar}:=
 \frac{1}{2}\bigl(
 \nabla_x+\sqrt{-1}\nabla_y
 \bigr)$
and 
$\del_{E,t}:=
 \nabla_t-\sqrt{-1}\phi$
on $E$.
The Bogomolny equation implies
the commutativity
$[\del_{E,\wbar},\del_{E,t}]=0$.
Thus, we obtain the mini-holomorphic structure
$(\del_{E,t},\del_{E,\wbar})$ on $E$.

\subsection{The induced instantons}
\label{subsection;16.12.16.11}

Let $\varphi:\cnum^2\lrarr\real\times\cnum$
be given by
$\varphi(u_1,u_2)=(|u_1|^2-|u_2|^2,2u_1u_2)$.
Let $U$ be any open subset in $\real\times\cnum$
such that $(0,0)\not\in U$.
Let $(E,h,\nabla,\phi)$ be any monopole on $U$.
We set
$(\Etilde,\htilde):=
\varphi^{\ast}(E,h)$
and 
\[
 \nablatilde:=
 \varphi^{\ast}\nabla
+\sqrt{-1}\xi\otimes\varphi^{\ast}\phi
\]
on $\varphi^{-1}(U)$,
where 
$\xi:=-u_1d\ubar_1+\ubar_1du_1-\ubar_2du_2+u_2d\ubar_2$.
As discovered by
Kronheimer \cite{Kronheimer-Master-Thesis}
(see also 
\cite{Charbonneau-Hurtubise},
\cite{Nash}
and \cite{Pauly}),
$(\Etilde,\htilde,\nablatilde)$
is an instanton on $\varphi^{-1}(U)$.
Namely,
the curvature $F(\nablatilde)$ is a $(1,1)$-form
and satisfies $\Lambda_{\cnum^2}F(\nablatilde)=0$.
This procedure induces an equivalence
between 
monopoles on $U$
and $S^1$-equivariant instantons
on $\varphi^{-1}(U)$.

We have the decomposition
$\nablatilde=\delbar_{\Etilde}\oplus\del_{\Etilde}$
into the $(0,1)$-part and the $(1,0)$-part,
and $\delbar_{\Etilde}$ gives a holomorphic structure
on $\Etilde$.
We can check the following lemma
by direct computations.
\begin{lem}
The holomorphic structure $\delbar_{\Etilde}$
is equivalent to the holomorphic structure
induced by the mini-holomorphic structure 
$(\del_{E,t},\del_{E,\wbar})$
as in {\rm\S\ref{subsection;16.12.16.10}}.
\hfill\qed
\end{lem}

\subsection{Dirac type singularity of monopoles}
\label{subsection;16.7.4.20}

Let $U\subset \real\times\cnum$ 
be an open subset such that
$(0,0)\in U$.
We set $U^{\ast}:=U\setminus\{(0,0)\}$.
Let $(E,\nabla,h,\phi)$ be a monopole on $U^{\ast}$.
We have the induced instanton
$(\Etilde,\htilde,\nablatilde)$
on $\varphi^{-1}(U^{\ast})$
as explained in \S\ref{subsection;16.12.16.11}.
It is standard to impose the following condition
on the behaviour of the monopole around the singularity $(0,0)$.

\begin{df}
\label{df;16.12.18.30}
If the induced instanton 
$(\Etilde,\htilde,\nablatilde)$ on $\varphi^{-1}(U^{\ast})$
is extended to an instanton
on $\varphi^{-1}(U)$,
then the singularity $(0,0)$ of $(E,\nabla,h,\phi)$
is called of Dirac type.
\hfill\qed
\end{df}

The concept of Dirac type singularity of monopoles 
was first introduced by Kronheimer
\cite{Kronheimer-Master-Thesis}.
See also 
\cite{Charbonneau-Hurtubise},
\cite{Cherkis-Kapustin2},
\cite{Nash} and \cite{Pauly}.
It is proved that 
the condition in Definition \ref{df;16.12.18.30}
is equivalent to some estimates for $\phi$
and its derivative.
In particular,
it is known that we have
$|\phi|_h=O\bigl((|t|+|w|)^{-1}\bigr)$
if $(0,0)$ is a Dirac type singularity of 
a monopole $(E,h,\nabla,\phi)$.
(See also Proposition \ref{prop;16.12.18.100} below.)
We shall prove the converse,
i.e., the estimate 
$|\phi|_h=O\bigl((|t|+|w|)^{-1}\bigr)$
implies that $(0,0)$ is a Dirac type singularity
of $(E,h,\nabla,\phi)$
(See Theorem \ref{thm;16.12.14.10}.)

\section{Characterizations of Dirac type singularity of monopoles}

\subsection{Characterization in terms of the growth order of 
 mini-holomorphic sections}
\label{subsection;16.12.14.1}

Let $U=U_t\times U_w\subset\real\times\cnum$
be an open subset such that
$(0,0)\in U$.
We set $U^{\ast}:=U\setminus\{(0,0)\}$.

Let $(E,h,\nabla,\phi)$ be a monopole of rank $r$
on $U^{\ast}$.
Let $(E,\del_{E,t},\del_{E,\wbar})$ be 
the underlying mini-holomorphic bundle
on $U^{\ast}$.
Take a small $\epsilon>0$.
Let $E^{\epsilon}$ and $E^{-\epsilon}$
denote the restriction of $E$
to $\{\epsilon\}\times U_w$
and $\{-\epsilon\}\times U_w$,
respectively.
Any holomorphic section $s^{-\epsilon}$
of $E^{-\epsilon}$ is naturally extended 
to a mini-holomorphic section $\stilde^{-\epsilon}$ of $E$
on $A^{-}:=U\setminus\{(t,0)\,|\,t\geq 0\}$.
Any holomorphic section $s^{\epsilon}$
of $E^{\epsilon}$ is naturally extended 
to a mini-holomorphic section $\stilde^{\epsilon}$ of $E$
on $A^{+}:=U\setminus\{(t,0)\,|\,t\leq 0\}$.
We put 
$U^+:=\{(t,w)\in U\,|\,t\geq 0\}$
and $U^-:=\{(t,w)\in U\,|\,t\leq 0\}$.
We consider the following condition {\bf (D)},
which is equivalent to the condition stated in
\S\ref{subsection;17.2.20.1}.
\begin{itemize}
\item
 For any holomorphic section $s^{\pm\epsilon}$
 of $E^{\pm\epsilon}$,
we have
 $\bigl|\stilde^{\pm\epsilon}\bigr|
\leq C(|t|^2+|w|^2)^{-N}$ 
 for some $C>0$  and $N>0$
 on $U^{\pm}\setminus\{(0,0)\}$.
\item
 The same estimate hold for sections of the dual of $(E,h,\nabla,\phi)$.
\end{itemize}

\begin{thm}
\label{thm;13.12.10.10}
The singularity $(0,0)$ of the monopole $(E,h,\nabla,\phi)$
is of Dirac type
if and only if the condition {\bf (D)} is satisfied.
\end{thm}
\pf
Suppose that 
the induced instanton $(\Etilde,\htilde,\nablatilde)$
on $\varphi^{-1}(U^{\ast})$
is extended to an instanton 
$(\Etilde_0,\htilde_0,\nablatilde_0)$
on $\varphi^{-1}(U)$.
By shrinking $U_w$,
we may assume to have the mini-holomorphic frames
$\sigma_{\pm,j}$ $(j=1,\ldots,r)$
of $E_{|A^{\pm}}$
as in the proof of Lemma \ref{lem;16.12.18.3}.
We have 
$\bigl| \sigma_{\pm,j}
 \bigr|_h
=O\bigl((|t|+|w|)^{-N}\bigr)$
for some $N>0$
on $U^{\pm}\setminus\{(0,0)\}$.
We also have similar frames and estimates
for the dual.
Then, we can easily check that 
the condition {\bf (D)} is satisfied.

Suppose the condition  {\bf (D)} is satisfied,
and we shall prove that 
$(0,0)$ is of Dirac type singularity of 
the monopole $(E,h,\nabla,\phi)$.
We have a holomorphic isomorphism
$\Phi:E^{-\epsilon}_{|U_w\setminus\{0\}}
\simeq
 E^{\epsilon}_{|U_w\setminus\{0\}}$.
\begin{lem}
\label{lem;16.12.19.20}
If the condition {\bf (D)} is satisfied, 
$\Phi$ is extended to a meromorphic isomorphism
$E^{-\epsilon}(\ast 0)\simeq E^{\epsilon}(\ast 0)$.
\end{lem}
\pf
For $\star=\pm$,
we take a holomorphic frames of
$v_1^{\star},\ldots,v_r^{\star}$
of $E^{\star\epsilon}$.
Let $(v_i^{\star})^{\lor}$ be the dual frames
of $(E^{\star\epsilon})^{\lor}$.
We have
$\Phi(v_i^{-})
=\sum_{j=1}^r
 \langle (v_j^{+})^{\lor},v_i^-\rangle
 v_j^{+}$.
Here,
$\langle\cdot,\cdot\rangle$
denote the pairing of 
sections of $E^{\lor}$ and $E$.
Note that
$\langle (v_j^{+})^{\lor},v_i^-\rangle$
are holomorphic on $U_t\times(U_w\setminus\{0\})$,
i.e., they are constant with respect to $t$,
and holomorphic with respect to $w$.
On $\{t=0\}$,
we have
$\Bigl|
 \langle (v_j^{+})^{\lor},v_i^-\rangle
 \Bigr|
=O(|w|^{-N})$
for some $N>0$.
Hence, we obtain that
$\langle (v_j^{+})^{\lor},v_i^-\rangle$
are meromorphic at $w=0$.
It means that
$\Phi$ is a meromorphic isomorphism.
\hfill\qed

\vspace{.1in}

Let $\varphi:\cnum^2\lrarr \real\times\cnum$
be given by
$\varphi(u_1,u_2)=\bigl(|u_1|^2-|u_2|^2,2u_1u_2\bigr)$.
Let us prove that 
the induced instanton
$(\Etilde,\htilde,\nablatilde)$
on $\varphi^{-1}(U^{\ast})$
is extended to an instanton on $\varphi^{-1}(U)$
under the condition {\bf (D)}.

By Lemma \ref{lem;16.12.19.20},
we have a holomorphic frame
$(s^-_1,\ldots,s^-_r)$ of $E^{-\epsilon}$
and a tuple of integers $\ell_1,\ldots,\ell_r$
such that
$(s^+_1,\ldots,s^+_r):=
 \bigl(
 w^{\ell_1}\Phi(s^-_1),\ldots,w^{\ell_r}\Phi(s^-_r)
 \bigr)$
is a holomorphic frame of $E^{\epsilon}$.
We have the following:
\begin{itemize}
\item
 We have 
 $|\stilde^{\pm}_i|\leq C(|t|^2+|w|^2)^{-N}$ 
 on $U^{\pm}\setminus\{(0,0)\}$
 for some $C>0$ and $N>0$.
\end{itemize}

On $\varphi^{-1}(U^{\ast})$,
we have a tuple of holomorphic sections
$e_i:=
 (2u_2)^{\ell_i}\varphi^{-1}(\stilde^-_i)
=u_1^{-\ell_i}\varphi^{-1}(\stilde^+_i)$
of $\Etilde$, which gives a frame.
Let $e_i^{\lor}$ denote the dual frame.
We have
$|e_i|_{\htilde}\leq C(|u_1|^2+|u_2|^2)^{-N}$
and 
$|e_i^{\lor}|_{\htilde^{\lor}}
 \leq C(|u_1|^2+|u_2|^2)^{-N}$
for some $C>0$ and $N>0$.
By the frame,
we extend $\Etilde$
to a holomorphic bundle $\Etilde_0$
on $\nbigu:=\varphi^{-1}(U)$.

Let us consider the case $\rank E=1$.
Note that
$\log \htilde(e_1,e_1)$ is a harmonic function
on $\nbigu\setminus\{(0,0)\}$
satisfying 
$\bigl|
 \log \htilde(e_1,e_1)
\bigr|
=O\Bigl(
 -\log(|u_1|^2+|u_2|^2)
 \Bigr)$.
The function $\log \htilde(e_1,e_1)$ is $L^2$,
and it is easy to check that
$\Delta \log \htilde(e_1,e_1)=0$ on $\nbigu$
as a distribution,
where
$\Delta:=
-(\del_{u_1}\del_{\ubar_1}+\del_{u_2}\del_{\ubar_2})$.
(See Lemma \ref{lem;16.12.18.2} below.)
Hence, we obtain that $\log \htilde(e_1,e_1)$
is a harmonic function on $\nbigu$
by the elliptic regularity.
In particular, it is $C^{\infty}$ on $\nbigu$.
Thus, the rank one case is proved.

Let us consider the general case.
Note that 
if the condition {\bf (D)} is satisfied for $(E,h,\nabla,\phi)$,
then the condition {\bf (D)} is also satisfied
for the determinant bundle of $(E,h,\nabla,\phi)$.
Hence, by applying the result
in the rank one case,
we obtain that
$\det(\htilde)$ gives a Hermitian-Einstein metric of
$\det(\Etilde_0)$ on $\nbigu$.
According to \cite[Theorem 1]{Donaldson-boundary-value},
we have a unique Hermitian-Einstein metric $\htilde_1$ of
$\Etilde_0$ on $\nbigu$
such that 
$\htilde_{1|\del\nbigu}
=\htilde_{|\del\nbigu}$.
By the uniqueness,
we have $\det \htilde_1=\det \htilde$.
Let $k$ be the endomorphism of $\Etilde$
determined by $\htilde=\htilde_1k$.
Note that $k$ is self-adjoint
with respect to both $\htilde$ and $\htilde_1$.
By using \cite[Lemma 3.1]{Simpson88},
we obtain that 
$\Delta\log \Tr(k)\leq 0$
on $\nbigu\setminus\{(0,0)\}$.
We also have
$\log\Tr(k)=O\bigl(
 -\log(|u_1|^2+|u_2|^2)
\bigr)$.
It is easy to check that
$\Delta\log\Tr(k)\leq 0$
as distributions on $\nbigu$.
(See Lemma \ref{lem;16.12.18.2} below.)
Namely, we obtain that 
$\log\Tr(k)$ 
is a subharmonic function on $\nbigu$.
Because we have
$\Tr(k)=\rank E$
on $\del\nbigu$,
we have
$\Tr(k)\leq \rank E$.
Because $\det (k)=1$
and because $k$ is a positive definite self-adjoint 
endomorphism of $(\Etilde,\htilde)$,
we obtain that $k=\id$.
Hence, $\htilde$ is $C^{\infty}$
on $\nbigu$,
and it is a Hermitian-Einstein metric 
of $\Etilde_0$.
\hfill\qed

\vspace{.1in}
We used the following lemma.
\begin{lem}
\label{lem;16.12.18.2}
Let $\nbigu$ be a neighbourhood of
$(0,0)$ in $\cnum^2$.
Let $f$ be an $\real$-valued $C^{2}$-function
on $\nbigu^{\ast}:=\nbigu\setminus\{(0,0)\}$
such that 
(i) $\Delta f\leq 0$ on $\nbigu^{\ast}$,
(ii) $|f|=O\bigl(\log(|u_1|^2+|u_2|^2)\bigr)$ around $(0,0)$.
Then, we have
$\Delta f\leq 0$ on $\nbigu$ as a distribution.
\end{lem}
\pf
We take a $C^{\infty}$-function
$\rho:\real\lrarr\real_{\geq 0}$ such that
$\rho(t)=0$ $(t\leq 1/2)$ and $\rho(t)=0$ $(t\geq 1)$.
For any  large $N$,
we put 
$\chi_N(u_1,u_2):=\rho\bigl(-N^{-1}\log(|u_1|^2+|u_2|^2)\bigr)$.
We have the following equalities:
\begin{equation}
 \del_{\ubar_1}\chi_N
=\rho'\bigl(-N^{-1}\log(|u_1|^2+|u_2|^2)\bigr)
 \cdot
 \frac{-u_1}{N(|u_1|^2+|u_2|^2)}
\end{equation}
\begin{multline}
  \del_{u_1}\del_{\ubar_1}\chi_N
=\rho''\bigl(-N^{-1}\log(|u_1|^2+|u_2|^2)\bigr)
 \cdot
 \frac{|u_1|^2}{N^2(|u_1|^2+|u_2|^2)^2}
 \\
-\rho'\bigl(-N^{-1}\log(|u_1|^2+|u_2|^2)\bigr)
 \cdot
 \frac{|u_2|^2}{N(|u_1|^2+|u_2|^2)^2}
\end{multline}
We have similar equalities for
$\del_{\ubar_2}\chi_N$
and $\del_{u_2}\del_{\ubar_2}\chi_N$.
Hence, we have a positive constant $C$
such that the following holds:
\begin{equation}
 \label{eq;16.12.18.1}
 \bigl|
 \del_{u_i}\chi_N
 \bigr|
\leq C(|u_1|+|u_2|)^{-1},
\quad\quad
 \bigl|
  \del_{\ubar_i}\del_{u_i}\chi_N
 \bigr|
\leq C(|u_1|+|u_2|)^{-2}.
\end{equation}

Let $\phi$ be any $\real_{\geq 0}$-valued test function on $\nbigu$.
We have 
$\int_{\nbigu}
 \Delta(f)\cdot \chi_N\phi\leq 0$ for any $N$.
We also have the following:
\begin{multline}
 \int_{\nbigu}
 \Delta(f)\cdot \chi_N\phi
=\int_{\nbigu}
 f\cdot \Delta(\chi_N\phi)
= \\
 \int_{\nbigu}
 f\cdot
 \Bigl( 
 (\Delta\chi_N)\cdot\phi
-\sum\del_{u_i}\chi_N\cdot\del_{\ubar_i}\phi
-\sum\del_{\ubar_i}\chi_N\cdot\del_{u_i}\phi
+\chi_N\Delta\phi
 \Bigr).
\end{multline}
By the estimate (\ref{eq;16.12.18.1})
and the assumption
$|f|=O\bigl(\log(|u_1|^2+|u_2|^2)\bigr)$,
we obtain 
$\int_{\nbigu}f\cdot\Delta\phi\leq 0$.
Hence, we obtain $\Delta(f)\leq 0$
on $\nbigu$.
\hfill\qed

\begin{example}
We need the estimate for the dual
in the condition {\bf (D)}.
Let $E$ be the product bundle
$\cnum\times U^{\ast}$
over $U^{\ast}$
with a prescribed frame $e$.
Set $R:=\sqrt{|t|^2+|w|^2}$.
Let $h$ be the Hermitian metric
$h(e,e)=\exp(-R^{-1})$.
We define $\nabla$ and $\phi$
by 
\[
 \nabla e
=e\cdot \Bigl(
 \frac{tdt}{2R^3}
+\frac{\wbar dw}{2R^3}
 \Bigr),
\quad\quad
 \phi=-\frac{\sqrt{-1}t}{2R^3}.
\]
Then, 
$(E,h,\nabla,\phi)$ is a monopole on $U^{\ast}$.
We have
$(\nabla_t-\sqrt{-1}\phi)e=0$
and $\nabla_{\wbar}e=0$,
i.e.,
$e$ is a mini-holomorphic frame of $E$.
We have 
$|e|_h=\exp(-R^{-1}/2)=O(1)$.
But, for the dual frame $e^{\lor}$,
we have
$|e^{\lor}|_{h^{\lor}}=\exp(R^{-1}/2)$
which is not dominated by
$(|t|^2+|w|^2)^{-N}$ for any $N$.

The induced holomorphic bundle 
$\Etilde$ on $\varphi^{-1}(U^{\ast})$
has a global holomorphic frame 
$\etilde:=\varphi^{-1}(e)$,
with which $\Etilde$ is extended to
a line bundle on $\varphi^{-1}(U)$.
We have 
$\htilde(\etilde,\etilde)=
 \exp\bigl(-(|u_1|^2+|u_2|^2)^{-1}\bigr)$.
Clearly, the metric is not $C^{\infty}$ on $\varphi^{-1}(U)$,
i.e.,
the induced instanton 
on $\varphi^{-1}(U^{\ast})$
is not extended across $(0,0)$.
\hfill\qed
\end{example}

\subsection{Characterization in terms of the growth order of
Higgs field}

Let $U$ and $U^{\ast}$ be as in \S\ref{subsection;16.12.14.1}.
Let $(E,h,\nabla,\phi)$ be a monopole on $U^{\ast}$.

\begin{thm}
\label{thm;16.12.14.10}
The singularity $(0,0)$ of 
the monopole $(E,h,\nabla,\phi)$
is of Dirac type
if and only if 
\[
|\phi|_h=O\bigl((|t|+|w|)^{-1}\bigr).
\]
\end{thm}
\pf
We use the notation in \S\ref{subsection;16.12.14.1}.
Take any holomorphic section $s^{-\epsilon}$
of $E^{-\epsilon}$,
and extend it to a mini-holomorphic section $\stilde^{-\epsilon}$
on $A^-$.
We have the following equality on $A^-$:
\[
 \frac{d}{dt}\bigl|\stilde^{-\epsilon}\bigr|_h^2
=2\Re h\bigl(
 \nabla_t\stilde^{-\epsilon},
 \stilde^{-\epsilon}
 \bigr)
=
 2\Re h\bigl(
 \sqrt{-1}\phi\stilde^{-\epsilon},
 \stilde^{-\epsilon}
 \bigr).
\]
By the assumption on $\phi$,
the following holds on $A^-$
for some $C_1>0$:
\[
 \Bigl|
 \frac{d}{dt}\bigl|\stilde^{-\epsilon}\bigr|_h^2
 \Bigr|
\leq
 \frac{C_1}{|t|+|w|}
 \bigl|\stilde^{-\epsilon}\bigr|_h^2.
\]
Hence, we obtain the following inequality on $A^-$:
\[
 \Bigl|
 \frac{d}{dt}\log|\stilde^{-\epsilon}|_h^2
 \Bigr|
\leq
 \frac{C_1}{|t|+|w|}.
\]
Thus, we obtain
$\bigl|
 \stilde^{-\epsilon}
 \bigr|_h
\leq
 C_2(|t|+|w|)^{-N}$ on 
$U^-\setminus\{(0,0)\}$
for some $C_2>0$ and $N>0$.

Similarly, for any holomorphic section
$s^{\epsilon}$ of $E^{\epsilon}$,
we extend it to
a mini-holomorphic section $\stilde^{\epsilon}$
of $E$ on $A^+$,
and then we have 
$\bigl|
 \stilde^{\epsilon}
 \bigr|_h\leq
 C_3(|t|+|w|)^{-N_3}$
on $U^+\setminus\{(0,0)\}$
for some $C_3>0$ and $N_3>0$.
We have similar estimates
for sections of the dual $(E,h,\nabla,\phi)^{\lor}$.
Then, the claim of Theorem \ref{thm;16.12.14.10}
follows from Theorem \ref{thm;13.12.10.10}.
\hfill\qed

\section{Asymptotic behaviour}

Although Dirac type singularity
is characterized by a rather weak condition 
as in Theorem \ref{thm;16.12.14.10},
monopoles are asymptotically close 
to the direct sum of Dirac monopoles
around their Dirac type singularity
\cite{Charbonneau-Hurtubise}.
We study it in a slightly more general situation.

\subsection{Hermitian metrics on mini-holomorphic bundles}
\label{subsection;16.12.18.110}

Let $(E,\del_{E,t},\del_{E,\wbar})$ be 
a mini-holomorphic bundle on 
any open subset $U\subset\real\times\cnum$.
Let $h$ be a Hermitian metric of $E$.
We have the unique differential operators
$\del'_{E,h,t}$ and $\del_{E,h,w}$ on $E$
satisfying the following conditions.
\begin{itemize}
\item
For any $f\in C^{\infty}(U)$
and $u\in C^{\infty}(E)$,
we have
$\del_{E,h,t}'(fu)=\del_t(f)\cdot u+f\del_{E,h,t}'u$
and 
$\del_{E,h,w}(fu)=
 \del_w(f)\cdot u+f\cdot\del_{E,h,w}u$.
\item
For any $u,v\in C^{\infty}(E)$,
we have
\[
 \del_t h(u,v)
=h(\del_{E,t}u,v)+h(u,\del'_{E,h,t}v),
\]
\[
 \del_{\wbar}h(u,v)
=h(\del_{E,\wbar}u,v)
+h(u,\del_{E,h,w}v).
\]
\end{itemize}
We set
$\nabla_{h,t}:=\frac{1}{2}(\del_{E,t}+\del'_{E,h,t})$,
$\nabla_{h,\wbar}:=\del_{E,\wbar}$
and $\nabla_{h,w}:=\del_{E,h,t}$.
We obtain the unitary connection
$\nabla_h$ on $(E,h)$
by $\nabla_h(s):=
 \nabla_{h,\wbar}(s)d\wbar
+\nabla_{h,w}(s)dw
+\nabla_{h,t}(s)dt$.
We also obtain the anti-self adjoint endomorphism
$\phi_h:=\frac{\sqrt{-1}}{2}(\del_{E,t}-\del'_{E,h,t})$ 
of $(E,h)$.

\vspace{.1in}

Suppose $(0,0)\not\in U$.
We have the induced holomorphic bundle
$(\Etilde,\delbar_{\Etilde})$
with the induced Hermitian metric $\htilde:=\varphi^{-1}(h)$.
We have the Chern connection $\nabla_{\htilde}$
of $(\Etilde,\delbar_{\Etilde},\htilde)$.
\begin{lem}
\label{lem;16.12.18.20}
For any $s\in C^{\infty}(E)$,
we have the following equalities:
\begin{equation}
 \label{eq;16.12.18.10}
  \varphi^{\ast}(\del_{E,h,t}'s)
=\frac{1}{|u_1|^2+|u_2|^2}
 \bigl(
 u_1\nabla_{\htilde,u_1}\varphi^{\ast}(s)
-u_2\nabla_{\htilde,u_2}\varphi^{\ast}(s)
 \bigr),
\end{equation}
\begin{equation}
 \label{eq;16.12.18.11}
 \varphi^{\ast}(\del_{E,h,w}s)
=\frac{1}{2(|u_1|^2+|u_2|^2)}
 \bigl(
 \ubar_2\nabla_{\htilde,u_1}\varphi^{\ast}(s)
+\ubar_1\nabla_{\htilde,u_2}\varphi^{\ast}(s)
 \bigr).
\end{equation}
\end{lem}
\pf
The right hand side of (\ref{eq;16.12.18.10}) and (\ref{eq;16.12.18.11})
are $S^1$-invariant.
By taking the descent,
we obtain the differential operators
$\del_{E,h,t}^{\circ}$ and $\del_{E,h,w}^{\circ}$
on $E$.
We can easily check that 
$\del_{E,h,t}^{\circ}$ and $\del_{E,h,w}^{\circ}$
satisfy the conditions for
$\del_{E,h,t}'$ and $\del_{E,h,w}$, respectively.
Then, the claim follows.
\hfill\qed

\subsection{Dirac monopoles}

Let $U$ be any neighbourhood of $(0,0)$
in $\real\times\cnum$.
We set $U^{\ast}:=U\setminus\{(0)\}$,
and 
$A^{\pm}:=U\setminus\{(t,0)\,|\,\pm t\leq 0\}$.
Let $L(m)$ be the mini-holomorphic bundle of rank one
on $U^{\ast}$
equipped with mini-holomorphic frames
$\sigma^{(m)}_{\pm}$
of $E_{|A^{\pm}}$
such that
$\sigma^{(m)}_-=(w/2)^{m}\sigma^{(m)}_+$
on $A^+\cap A^-$.
The induced holomorphic line bundle
$\widetilde{L(m)}_0$ on $\varphi^{-1}(U)$
is equipped with the frame $e^{(m)}$
such that $b^{\ast}(e^{(m)})=b^me^{(m)}$ for any $b\in S^1$,
and that $u_1^{-m}e^{(m)}$ and $u_2^{m}e^{(m)}$
induce $\sigma^{(m)}_{+}$
and $\sigma_-^{(m)}$,
respectively.

Let $\htilde_0^{(m)}$ be the metric of $\widetilde{L(m)}_0$
given by $\htilde_0^{(m)}(e^{(m)},e^{(m)})=1$.
It induces a metric $h^{(m)}$ of $L^{(m)}$.
We obtain the unitary connection
$\nabla^{(m)}$
and the anti-self adjoint endomorphism
$\phi^{(m)}$
by the procedure in \S\ref{subsection;16.12.18.110}.
By the construction,
the tuple $(L^{(m)},h^{(m)},\nabla^{(m)},\phi^{(m)})$
is a monopole,
called the Dirac monopole.
It is easy to check that
$\phi^{(m)}$ is the multiplication of
$\frac{\sqrt{-1}m}{R}$,
where $R=\sqrt{|w|^2+|t|^2}$.

\subsection{Extendability}

Let $U$ be a neighbourhood of $(0,0)\in\real\times\cnum$.
We set $U^{\ast}:=U\setminus\{(0,0)\}$.
Let $(E,\del_{E,t},\del_{E,\wbar})$
be a mini-holomorphic bundle such that
the induced holomorphic bundle
$(\Etilde,\delbar_{\Etilde})$ on
$\varphi^{-1}(U^{\ast})$
is extended to a holomorphic bundle
$(\Etilde_0,\delbar_{\Etilde_0})$
on $\varphi^{-1}(U)$.
We set $R:=(|t|^2+|w|^2)^{1/2}$.

Let $h$ be a Hermitian metric of $E$
such that $\varphi^{-1}(h)$
induces a $C^{\infty}$-metric $\htilde_0$
on $\Etilde_0$.
We have the unitary connection $\nabla_h$
and the anti-self adjoint endomorphism $\phi_h$
of $E$ as in \S\ref{subsection;16.12.18.110}.

\begin{prop}
\label{prop;16.12.18.100}
We have an isomorphism of 
mini-holomorphic bundles
$\Phi:(E,\del_{E,t},\del_{E,\wbar})
\simeq
 \bigoplus_{i=1}^{\rank E}
 L(k_i)$
such that the following holds.
\begin{description}
\item[(P1)]
 Set $h_1:=\bigoplus h^{(k_i)}$.
 We have the endomorphism $a$ of $E$
 determined by 
 $h=\Phi^{\ast}(h_1)\cdot a$.
 Then, $|a-\id|_h=O(R)$.
\item[(P2)]
 We set $\phi_1:=\bigoplus \phi^{(k_i)}$.
 Then,
 $\phi_h-\Phi^{\ast}(\phi_1)$
 is bounded with respect to $h$.
\item[(P3)]
 We set $\nabla_1:=\bigoplus \nabla^{(k_i)}$.
 Then,
 $\nabla_h-\Phi^{\ast}(\nabla_1)$
 is bounded with respect to $h$.
\end{description}
Moreover,
$\nabla_h(R\phi_h)$ is bounded
with respect to $h$.
\end{prop}
\pf
We begin with the study of
holomorphic frames of $\Etilde_0$.

\begin{lem}
\label{lem;16.12.18.10}
We have a holomorphic frame 
$e_1,\ldots,e_{\rank E}$ of $\Etilde_0$
satisfying the following conditions:
\begin{itemize}
\item
 $b^{\ast}e_i=b^{k_i}e_i$
 for some integers $k_i$
 for any $b\in S^1$.
\item
 $\htilde_0(e_i,e_j)-\delta_{i,j}=O(|u_1|^2+|u_2|^2)$,
 where $\delta_{i,i}=1$ and $\delta_{i,j}=0$ $(i\neq j)$.
\end{itemize}
\end{lem}
\pf
We can take a holomorphic frame 
$e'_1,\ldots,e'_{\rank E}$ satisfying the first condition
by the argument in Lemma \ref{lem;16.12.18.3}.
We also assume that 
$(e'_1,\ldots,e'_{\rank E})_{|(0,0)}$ is an orthonormal  frame of
$(\Etilde_0,\htilde_0)_{|(0,0)}$.
Let us observe that
we can modify it so that the second condition is also satisfied.

We have 
$b^{\ast}\bigl(
 \htilde_0(e'_i,e'_j)\bigr)=b^{k_i-k_j}\htilde_0(e'_i,e'_j)$
for any $b\in S^1$.
We have the Taylor expansion of $\htilde_0(e_i,e_j)$:
\[
 \sum_{a,b,c,d\geq 0}
 G_{i,j;a,b,c,d}\cdot u_1^a\ubar_1^bu_2^c\ubar_2^d.
\]
We have
$G_{i,j;a,b,c,d}=0$ unless
$a-b-c+d=k_i-k_j$.
We also have
$G_{i,j;a,b,c,d}
 =\overline{G_{j,i;b,a,d,c}}$.

Suppose that 
$a-b-c+d=k_i-k_j$.
If $|k_i-k_j|\geq 2$,
we have
$|a|+|b|+|c|+|d|\geq 2$.
If $k_i-k_j=0$,
we have 
$(a,b,c,d)=(0,0,0,0)$
or 
$|a|+|b|+|c|+|d|\geq 2$.
If $k_i-k_j=1$,
we have
$(a,b,c,d)=(1,0,0,0),(0,0,0,1)$
or $|a|+|b|+|c|+|d|\geq 2$.
If $k_i-k_j=-1$,
we have
$(a,b,c,d)=(0,1,0,0),(0,0,1,0)$
or $|a|+|b|+|c|+|d|\geq 2$.

We set
\[
 e_i:=e'_i
-\sum_{\{j;k_i-k_j=1\}}
 G_{i,j;1,0,0,0}\cdot u_1\cdot e'_j
-\sum_{\{j;k_i-k_j=-1\}}
 G_{i,j;0,0,1,0}\cdot u_2\cdot e'_j.
\]
Then, we can easily check that 
$e_1,\ldots,e_{\rank E}$ is a holomorphic frame
with the desired property.
\hfill\qed

\vspace{.1in}

Let $\vece=(e_1,\ldots,e_{\rank E})$ 
be a frame as in Lemma \ref{lem;16.12.18.10}.
From the decomposition
$\Etilde_0=\bigoplus\nbigo e_i$,
we obtain an isomorphism 
of mini-holomorphic bundles
$\Phi:E\simeq \bigoplus L(k_i)$.
We shall prove that 
$\Phi$ has the desired property.
It satisfies {\bf(P1)} by the above construction.

We have the Chern connection 
$\nabla_{\htilde_0}$
of $(\Etilde_0,\htilde_0)$.
Let $C_i$ $(i=1,2)$ be the matrix valued functions
determined by
$\nabla_{\htilde_0,u_i}\vece
=\vece\cdot C_i$.
We have $C_{i|(0,0)}=0$.

On $\varphi^{-1}(A^+)$,
we consider the frame
$v^+_{i}:=u_1^{-k_i}e_i$ $(i=1,\ldots,\rank E)$.
Because they are $S^1$-invariant,
they induce a frame 
$\sigma^+_i$ $(i=1,\ldots,\rank E)$.
We remark the following,
which is clear by the construction of $\vecsigma^+$.
\begin{lem}
\label{lem;16.12.18.200}
Let $G$ be an endomorphism of $E_{|A^+}$,
and let $B$ be the matrix valued function
determined by
$G\vecsigma^+=\vecsigma^+B$,
i.e.,
$G(\sigma^+_j)=\sum B_{i,j}\sigma^+_i$.
Then, 
$G$ is bounded with respect to $h$
if and only if 
$|B_{i,j}u_1^{k_j-k_i}|$
are bounded.
\hfill\qed
\end{lem}

Let $C^+_i$ $(i=1,2)$
be the matrix valued functions determined by
$C^+_{i;p,q}:=C_{i;p,q}u_1^{k_p-k_q}$.
Note that 
$u_iC^+_i$ $(i=1,2)$,
 $\ubar_2C^+_1$ and $\ubar_1C^+_2$
are $S^1$-invariant.
Hence, we regard them as 
matrix valued functions on $A^+$.
Let $\Gamma$ denote the diagonal matrix
whose $(i,i)$-components are $-k_i$.
We have
$\del_{E,t}\vecsigma^+=0$
and
$\del_{E,\wbar}\vecsigma^+=0$.
By using Lemma \ref{lem;16.12.18.20},
we have 
\[
 \del_{E,t}'\vecsigma^+
=\vecsigma^+
 \frac{1}{R}\bigl(u_1C_1^++\Gamma-u_2C_2^+\bigr),
\quad
 \del_{E,w}\vecsigma^+
=\vecsigma^+
 \frac{1}{2R}
 \bigl(\ubar_2C_1^++\ubar_2u_1^{-1}\Gamma+\ubar_1C_2^+\bigr).
\]
We set
$D_1:=u_1C_1^+-u_2C_2^+$
and $D_2:=\ubar_2C_1^++\ubar_1C_2^+$.
Then, we have
\[
 \phi_h\vecsigma^+=\vecsigma^+\cdot 
 \frac{\sqrt{-1}}{2R}(-D_1-\Gamma),
\quad
 \nabla_t\vecsigma^+=\vecsigma^+\frac{1}{2R}(D_1+\Gamma),
\]
\[
  \nabla_{\wbar}\vecsigma^+=0,
\quad
 \nabla_w\vecsigma^+=
 \vecsigma^+\cdot\frac{1}{2R}
 \bigl(D_2+\ubar_2u_1^{-1}\Gamma\bigr).
\]
We also have the following:
\[
 \Phi^{\ast}(\phi_1)\vecsigma^+
=\vecsigma^+\cdot
 \frac{\sqrt{-1}}{2R}(-\Gamma),
\quad
 \Phi^{\ast}(\nabla_{1,t})\vecsigma^+
=\vecsigma^+\frac{1}{2R}\Gamma,
\]
\[
  \Phi^{\ast}(\nabla_{1,\wbar})\vecsigma^+=0,
\quad
 \Phi^{\ast}(\nabla_{1,w})\vecsigma^+=
  \vecsigma^+\cdot\frac{1}{2R}
 \ubar_2u_1^{-1}\Gamma.
\]
By using Lemma \ref{lem;16.12.18.200},
we obtain the boundedness of
$\phi_h-\Phi^{\ast}(\phi_{1})$ and
$\nabla_h-\Phi^{\ast}(\nabla_{1})$
on $A^+$.

Because 
$\ubar_2\del_{u_1}D_1$ and
$\ubar_1\del_{u_2}D_1$ are $S^1$-invariant,
we may naturally regard them 
as matrix valued functions on $A^+$.
We also have
$\bigl(\ubar_2\del_{u_1}D_1\bigr)_{p,q}u^{k_q-k_p}
=O(R)$
and
$\bigl(\ubar_1\del_{u_2}D_1\bigr)_{p,q}u^{k_q-k_p}
=O(R)$.
We have the following formula:
\[
 \nabla_{w}(R\phi_h)\vecsigma^+
=\vecsigma^+
 \cdot
 \Bigl(
-\frac{\sqrt{-1}}{4R}
 [D_2+\ubar_2u_1^{-1}\Gamma,D_1+\Gamma]
-\frac{\sqrt{-1}}{2R}
 \bigl(
 \ubar_2\del_{u_1}D_1
+\ubar_1\del_{u_2}D_1
 \bigr)
 \Bigr).
\]
Note $[\ubar_2u_1^{-1}\Gamma,\Gamma]=0$.
Hence, we obtain that
$\nabla_{w}(R\phi_h)$ is bounded
on $A^+$.
Similarly, we obtain the boundedness of
$\nabla_t(R\phi_h)$
and $\nabla_{\wbar}(R\phi_h)$
on $A^+$.
Hence,
$\nabla(R\phi_h)$ is bounded on $A^+$.
We can obtain the estimate on $A^-$
in a similar way.
\hfill\qed

\paragraph{Acknowledgement}

The authors thank Sergey Cherkis 
for answering to our question about 
their characterization of Dirac type singularity 
in \cite{Cherkis-Kapustin2}.
The first author is partially supported by
the Grant-in-Aid for Scientific Research (S) (No. 24224001),
the Grant-in-Aid for Scientific Research (S) (No. 16H06335),
and the Grant-in-Aid for Scientific Research (C) (No. 15K04843),
Japan Society for the Promotion of Science.

\noindent
\begin{tabular}{ll}
(Takuro Mochizuki) &
 Research Institute for Mathematical Sciences,
 Kyoto University,\\
 &
Kyoto 606-8502, Japan, 
E-mail address:  takuro@kurims.kyoto-u.ac.jp\\
\mbox{{}}\\
(Masaki Yoshino) &
 Research Institute for Mathematical Sciences,
 Kyoto University, \\
 &Kyoto 606-8502, Japan,
E-mail address:  yoshino@kurims.kyoto-u.ac.jp
\end{tabular}

\end{document}